\documentclass[dvips, 12pt, a3paper]{article}
\usepackage{latexsym,amsmath,amssymb,graphicx}
\newtheorem{definition}{Definition}[section]
\newtheorem{theorem}[definition]{Theorem}

\newtheorem{lemma}[definition]{Lemma}
\begin{document}
\begin{flushleft}
 {\bf\Large {Characterization of Nonhomogeneous Wavelet Bi-frames for Reducing Subspaces of Sobolev Spaces over Local Fields of Prime Characteristics }}

\parindent=0mm \vspace{.1in}

\end{flushleft}

\vspace{.2in}\parindent=0mm
\begin{center}
  {{M. YOUNUS BHAT}}

\parindent=0mm \vspace{.1in}
{{\it Department of  Mathematical Sciences,  Islamic University of Science and Technology Awantipora, Pulwama, Jammu and Kashmir 192122, India. E-mail: $\text{gyounusg@gmail.com}$}}
\end{center}

\parindent=0mm \vspace{.1in}
{\bf{Abstract:}}  In this paper, we provide the characterization of nonhomogeneous wavelet bi-frames. First of all we introduce the reducing subspaces of Sobolev spaces over local fields of prime characteristics and then characterize the nonhomogeneous wavelet bi-frames over such fields.
\parindent=0mm \vspace{.1in}
{\bf{Keywords:}} Nonhomogeneous, Local field, Wavelet Bi-frames, Sobolev Spaces.

\parindent=0mm \vspace{.1in}
{\bf{Mathematics Subject Classification:}} 42C40; 42C15; 43A70; 11S85; 47A25.
\section{Introduction}
 
\parindent=8mm \vspace{.1in}
During the last two decades, there is a substantial body of work that has been concerned with the construction of wavelets on  local fields.  Even though the structures and metrics of local fields of zero and positive characteristics are similar, their wavelet and MRA (multiresolution analysis) theory are quite different.  The theory of wavelets, wavelet frames,  multiresolution analysis, Gabor frames on local fields of positive characteristics (LFPC) are extensively studied by many researchers including  Benedetto, Behera and Jahan,  Ahmad and Shah, Jiang, Li and Ji, Shukla and Mittal   in the references \cite{{ow2020},{nuwf},{ogbr},{BJ1},{BJ2},{BB},{oJGP},{oSIS},{SAJ},{oATA},{shukla2019},{shukla2018}} but still more concepts required to be studied for its enhancement on LFPC. Continuing our study of frames on local fields, we introduce a  comprehensive theory of discrete wave packet systems on local fields  by exploiting the machinery of  Fourier transforms and we also give the definition of discrete periodic wave packet transform. A characterization of the system to be a Parseval frame and discrete periodic wave packet frame for $\ell^2(\mathcal Z)$ are obtained. Furthermore,  we construct discrete periodic wave packet frame by filter sequence iteration. 

\parindent=8mm \vspace{.1in}
We in this paper,  provide the characterization of nonhomogeneous wavelet bi-frames. First of all we introduce the reducing subspaces of Sobolev spaces over local fields of prime characteristics and then characterize the nonhomogeneous wavelet bi-frames over such fields.

\parindent=8mm \vspace{.1in}
The rest of the paper is tailored as follows. In Section 2, we recall some basic Fourier analysis on local  fields  and also some results which are required in the subsequent sections. In section 3, Characterization of NWBFs over Local Field. In Section 4, we construct the $J^{th}$-stage discrete periodic wave packet frame for $\ell^2(\mathcal Z)$. 

\section{Preliminaries on Local Fields}
 \subsection{Local Fields}
A  local  field $K$ is a  locally compact, non-discrete and totally disconnected field. If it is of characteristic zero, then  it is a field of $p$-adic numbers $\mathbb Q_p$ or its finite extension. If $K$ is of positive characteristic, then $K$ is a field of formal Laurent series over a finite field $GF(p^c)$. If $c =1$, it is a $p$-series field, while for $c\ne 1$, it is an algebraic extension of degree  $c$ of a $p$-series field. Let $K$ be a fixed local  field  with the ring of integers ${\mathfrak D}= \left\{x \in K: |x| \le 1\right\}$. Since $K^{+}$ is a locally compact Abelian group, we choose a Haar measure $dx$ for $K^{+}$. The  field $K$ is locally compact, non-trivial, totally disconnected and complete topological field endowed with non--Archimedean norm  $|\cdot|:K\to \mathbb R^+$ satisfying

\parindent=0mm \vspace{.1in}
(a) $|x|=0$ if and only if $x = 0;$

\parindent=0mm \vspace{.1in}
(b) $|x\,y|=|x||y|$ for all $x, y\in K$;

\parindent=0mm \vspace{.1in}
(c) $|x+y|\le \max \left\{ |x|, |y|\right\}$ for all $x, y\in K$.

\parindent=0mm \vspace{.1in}
Property (c) is called the ultrametric inequality. Let ${\mathfrak B}= \left\{x \in K: |x| < 1\right\}$ be the prime ideal of the ring of integers ${\mathfrak D}$ in $K$. Then, the residue space ${\mathfrak D}/{\mathfrak B}$ is isomorphic to a finite field $GF(q)$, where $q = p^{c}$ for some prime $p$ and $c\in\mathbb N$. Since  $K$ is totally disconnected and $\mathfrak B$ is both prime and principal ideal, so there exist a prime element $\mathfrak p$ of $K$ such that ${\mathfrak B}= \langle \mathfrak p \rangle=\mathfrak p {\mathfrak D}$. Let ${\mathfrak D}^*= {\mathfrak D}\setminus {\mathfrak B }=\left\{x\in K: |x|=1   \right\}$. Clearly,  ${\mathfrak D}^*$ is a group of units in $K^*$ and if $x\not=0$, then can write $x=\mathfrak p^n y, y\in {\mathfrak D}^*.$ Moreover, if ${\cal U}= \left\{a_m:m=0,1,\dots,q-1 \right\}$ denotes the fixed full set of coset representatives of ${\mathfrak B}$ in ${\mathfrak D}$, then every element $x\in K$ can be expressed uniquely  as $x=\sum_{\ell=k}^{\infty} c_\ell \,\mathfrak p^\ell $ with $c_\ell \in {\cal U}.$ Recall that ${\mathfrak B}$ is compact and open, so each  fractional ideal ${\mathfrak B}^k= \mathfrak p^k {\mathfrak D}=\left\{x \in K: |x| < q^{-k}\right\}$  is also compact and open and is a subgroup of $K^+$. We use the notation in Taibleson's book \cite{tab}. In the rest of this paper, we use the symbols $\mathbb N, \mathbb N_0$ and $\mathbb Z$ to denote the sets of natural, non-negative integers and integers, respectively.

\parindent=8mm \vspace{.1in}
 Let $\chi$ be a fixed character on $K^+$ that is trivial on ${\mathfrak D}$ but  non-trivial on  ${\mathfrak B}^{-1}$. Therefore, $\chi$ is constant on cosets of ${\mathfrak D}$ so if $y \in {\mathfrak B}^k$, then $\chi_y(x)=\chi(y,x), x\in K.$ Suppose that $\chi_u$ is any character on $K^+$, then the restriction $\chi_u|{\mathfrak D}$ is a character on ${\mathfrak D}$. Moreover, as characters on ${\mathfrak D}, \chi_u=\chi_v$ if and only if $u-v\in {\mathfrak D}$. Hence, if  $\left\{u(n): n\in\mathbb N_0\right\}$ is a complete list of distinct coset representative of ${\mathfrak D}$ in $K^+$, then, as it was proved in [13], the set  $\left\{\chi_{u(n)}: n\in\mathbb N_0\right\}$   of distinct characters on ${\mathfrak D}$ is a complete orthonormal system on ${\mathfrak D}$.

\parindent=8mm \vspace{.1in}
We now impose a natural order on the sequence $\{u(n)\}_{n=0}^\infty$. We have ${\mathfrak D}/ \mathfrak B \cong GF(q) $ where $GF(q)$ is a $c$-dimensional vector space over the field $GF(p)$. We choose a set $\left\{1=\zeta_0,\zeta_1,\zeta_2,\dots,\zeta_{c-1}\right\}\subset {\mathfrak D^*}$ such that span$\left\{\zeta_j\right\}_{j=0}^{c-1}\cong GF(q)$. For $n \in \mathbb N_0$ satisfying
$$0\leq n<q,~~n=a_0+a_1p+\dots+a_{c-1}p^{c-1},~~0\leq a_k<p,~~\text{and}~k=0,1,\dots,c-1,$$

\parindent=0mm \vspace{.1in}
we define
$$u(n)=\left(a_0+a_1\zeta_1+\dots+a_{c-1}\zeta_{c-1}\right){\mathfrak p}^{-1}.$$

\parindent=0mm \vspace{.1in}
Also, for $n=b_0+b_1q+b_2q^2+\dots+b_sq^s, ~n\in \mathbb N_{0},~0\leq b_k<q,k=0,1,2,\dots,s$, we set

$$u(n)=u(b_0)+u(b_1){\mathfrak p}^{-1}+\dots+u(b_s){\mathfrak p}^{-s}.$$

\parindent=0mm \vspace{.1in}
This defines $u(n)$ for all $n\in \mathbb N_{0}$. In general, it is not true that $u(m + n)=u(m)+u(n)$. But, if $r,k\in\mathbb N_{0}\; \text{and}\;0\le s<q^k$, then $u(rq^k+s)=u(r){\mathfrak p}^{-k}+u(s).$ Further, it is also easy to verify that $u(n)=0$ if and only if $n=0$ and $\{u(\ell)+u(k):k \in \mathbb N_0\}=\{u(k):k \in \mathbb N_0\}$ for a fixed $\ell \in \mathbb N_0.$ Hereafter we use the notation $\chi_n=\chi_{u(n)}, \, n\ge 0$.

\parindent=8mm \vspace{.2in}
Let the local field $K$ be of characteristic $p>0$ and $\zeta_0,\zeta_1,\zeta_2,\dots,\zeta_{c-1}$ be as above. We define a character $\chi$ on $K$ as follows:
$$\chi(\zeta_\mu {\mathfrak p}^{-j})= \left\{
\begin{array}{lcl}
\exp(2\pi i/p),&&\mu=0\;\text{and}\;j=1,\\
1,&&\mu=1,\dots,c-1\;\text{or}\;j \neq 1.
\end{array}
\right. $$

\parindent=0mm \vspace{.1in}
\subsection{ Fourier Transforms on Local Fields}

\parindent=8mm \vspace{.0in}
 The Fourier transform of $f \in L^1(K)$ is denoted by $\widehat f(\omega)$ and defined  by
\begin{align*}
{F}\big\{f(x)\big\}=\widehat f(\omega)=\int_K f(x)\overline{ \chi_\omega(x)}\,dx.
\end{align*}
It is noted that
$$\widehat f(\omega)= \displaystyle \int_K f(x)\,\overline{ \chi_\omega(x)}dx= \displaystyle \int_K f(x)\chi(-\omega x)\,dx.$$

\parindent=8mm \vspace{.0in}
The properties of Fourier transforms on local field $K$ are much similar to those of on the classical field $\mathbb R$. In fact, the Fourier transform on local fields of positive characteristic have the following properties:
\begin{itemize}
  \item The map $f\to \widehat f$ is a bounded linear transformation of $L^1(K)$ into $L^\infty(K)$, and $\big\|\widehat f\big\|_{\infty}\le \big\|f\big\|_{1}$.
  \item If $f\in L^1(K)$, then $\widehat f$ is uniformly continuous.
  \item If $f\in L^1(K)\cap L^2(K)$, then $\big\|\widehat f\big\|_{2}=\big\|f\big\|_{2}$.
\end{itemize}

\parindent=8mm \vspace{.0in}
The Fourier transform of a function $f\in L^2(K)$ is defined by
\begin{align*}
\widehat f(\omega)= \lim_{k\to \infty} \widehat f_{k}(\omega)=\lim_{k\to \infty}\int_{|x|\le q^{k}} f(x)\overline{ \chi_\omega(x)}\,dx,
\end{align*}

\parindent=0mm \vspace{.0in}
where  $f_{k}=f\,\Phi_{-k}$ and $\Phi_{k}$ is  the characteristic function of ${\mathfrak B}^{k}$.  Furthermore, if $f\in L^2(\mathfrak D)$, then we define the Fourier coefficients of $f$ as
\begin{align*}
\widehat f\big(u(n)\big)=\int_{\mathfrak D} f(x) \overline{ \chi_{u(n)}(x)}\,dx.
\end{align*}

\parindent=8mm \vspace{.0in}
The series $\sum_{n\in \mathbb N_{0}} \widehat f\big( u(n)\big) \chi_{u(n)}(x)$ is called the Fourier series of $f$. From the standard $L^2$-theory for compact Abelian groups, we conclude that the Fourier series of $f$ converges to $f$ in $L^2(\mathfrak D)$ and Parseval's identity holds:
\begin{align*}
\big\|f\big\|^2_{2}=\int_{\mathfrak D}\big|f(x)\big|^2 dx= \sum_{n\in \mathbb N_{0}} \left| \widehat f\big(u(n)\big)\right|^2.
\end{align*}

\parindent=8mm \vspace{.0in}
Let $\mathcal Z =\{u(n) : n \in \mathbb N_0\},$ where $\{u(n) : n \in \mathbb N_0\}$ is a complete list of (distinct) coset representation of $\mathfrak D$ in $K^+.$ Then we define

\parindent=0mm \vspace{.0in}
$$\ell^2(\mathcal Z) = \left\{ z: {\mathcal Z}\rightarrow \mathbb C :\sum_{n \in \mathbb N_0}|z(u(n))|^2 < \infty\right\}$$

is a Hilbert space with inner product

\parindent=0mm \vspace{.0in}
$$\left\langle z,  w \right\rangle =\sum_{n \in \mathbb N_0} z(u(n))\overline{w(u(n))}.$$
The following definitions are natural which are used in the subsquent sections:
\begin{definition} The {\it Fourier transform} on $\ell^2(\mathcal Z)$ is the map  $ \widehat{} : \ell^2(\mathcal Z) \rightarrow L^2(\mathfrak D)$ defined for $z \in \ell^2(\mathcal Z)$ by

$$\widehat{z}(\omega) =\sum_{n \in \mathbb N_0} z(u(n))\chi_{u(n)}(\omega),$$
 and the {\it Inverse Fourier transform} on $L^2(\mathfrak D)$ is the map $\vee :L^2(\mathfrak D) \rightarrow \ell^2(\mathcal Z)$ defined for $ f \in L^2(\mathfrak D)$ by

$$f^{\vee}(u(n)) = \left\langle f , \chi_{u(n)}\right\rangle = \int_{\mathfrak D} f(x) \overline{\chi_{u(n)}(x)}\,dx.$$
\end{definition}

\parindent=8mm \vspace{.0in}
For $ z \in \ell^2(\mathcal Z),$ we have 
\begin{align*}
\left(\widehat z\right)^{\vee}(u(n)) & = \left\langle \widehat z, \chi_{u(n)}\right\rangle \\ 
& =\left \langle \displaystyle\sum_{m \in \mathbb N_0} z(u(m))\chi_{u(m)},\chi_{u(n)}\right\rangle\\ 
& = \displaystyle\sum_{m \in \mathbb N_0} z(u(m)) \left \langle \chi_{u(m)},\chi_{u(n)}\right\rangle\\
& =z(u(n)).
\end{align*}
Since $\left\{ \chi_{u(n)} : n \in \mathbb N_0\right\}$ is an orthonormal basis for $L^2(\mathfrak D)$. It is also clear that the function $\widehat{z}$ is an integral periodic function because for $ m \in \mathbb N_0,$ we have
\begin{align*}
\widehat{z}(\omega +u(m))&=\displaystyle\sum_{n \in \mathbb N_0} z(u(n))\chi_{u(n)}(\omega),\chi_{u(n)}(u(m)).\\
& = \displaystyle\sum_{n \in \mathbb N_0} z(u(n))\chi_{u(n)}(\omega)\\
& = \widehat{z}(\omega).
\end{align*}

\parindent=8mm \vspace{.0in}
For $z,w \in \ell^2(\mathcal Z),$ we have {\it Parseval's relation}:
$$\left\langle z, w\right\rangle =\sum_{n \in \mathbb N_0} z(u(n))\overline{w(u(n))} = \int_{\mathfrak{D}}\widehat{z}(\omega)\overline{\widehat{w}(\omega)}\,d\omega =\left\langle \widehat z,\widehat w\right\rangle,$$
and {\it Plancherel's formula}:
$$\|z\|^2 = \sum_{n \in \mathbb N_0} |z(u(n))|^2 = \int_{\mathfrak{D}}|\widehat{z}(\omega)|^2\,d\omega = \|\widehat{z}\|^2.$$
\begin{definition} \label{jjd2.1}
Let $\mathcal  H$ denotes a separable Hilbert space and $\mathcal{P}$ be countable index set. A countable set $\left\{\varphi_{_{P}} :p\in \mathcal{P}\right\}$ in $\mathcal{H}$ is called a frame $\mathcal{H}$ if there exists constants $A$ and $B$,
$0<C \le D< \infty$ such that
\begin{equation}\label{jjdee2.4} 
C\Vert f \Vert ^2 \le \sum\limits_{p \in \mathcal{P}} \vert \langle f, \varphi_{p} \rangle\vert^2\le D \Vert f \Vert^2
\end{equation}
holds for all $f \in \mathcal H$. The greatest possible such $C$ is the lower frame bound and the least possible such $D$ is the upper frame bound. We say that $\left\{\varphi_p :p \in \mathcal{P}\right\}$ is a tight frame for $\mathcal{H}$ with frame bound $C$ if $C = D$, and a normalized tight (Parseval) frame
if $C = D = 1$.
\end{definition}
\begin{definition}\label{jjd2.2}
For $m, k \in \mathbb N_0,$ the {\it translation operator} $T_{u(m)} :\ell^2(\mathcal{Z}) \to \ell^2(\mathcal{Z})$ is defined by 
$$T_{u(m)} z(u(n))=z\left(u(n)-u(m)\right)$$
and the {\it modulation operator} $M_{u(k)}:\ell^2(\mathcal{Z}) \to \ell^2(\mathcal{Z})$ by
 $$M_{u(k)}z\left(u(n)\right)=z\left(u(n)\right)\overline{\chi_{u(k)}}.$$
\end{definition}

\begin{definition}\label{jjd2.4}
 For $v_\ell ^j \in \ell^2(\mathcal{Z}), ~0 \le \ell \le L, 0 \le j\le J$, we call $\mathcal{W}(V)$ a $J^{th}$ -stage wave packet system associated with $ V=\left\{ v_{\ell}: 0 \le \ell \le L \right\}$ if 
\begin{equation}\label{jjd2.5}
\mathcal{W}(V)=\left\{\mathcal{W}_{k,m}^{j}v_{\ell}^{j}=T_{\mathfrak{p}^ju(m)}M_{u(k)}v^{j}_{\ell}: k \in \mathbb N_0, m \in \mathbb N_0, 0\le \ell \le L, 0\le j \le J\right\}.
\end{equation} 
 
\end{definition}

\parindent=0mm \vspace{.0in}
 For $J=0$, the wave packet system in (\ref{jjd2.5}) becomes  the Gabor system \\
$$\left\{ \overline{\chi_{u(k)}(u(m)} v_{\ell}(\cdot- u(m)): k \in  \in \mathbb N_0, 0 \le \ell \le L\right\}.$$

\parindent=0mm \vspace{.0in}
On taking $k=0$ in (\ref{jjd2.5}), the wave packet system degenerates into the wavelet system
$$\left\{v_{\ell}^{j}(\cdot- \mathfrak{p}^j u(m)): m  \in \mathbb N_0, 0\le \ell \le L, j=1,2,...,J \right\}.$$
Thus the  Gabor system and wavelet system are the particular cases of wave packet systems.
\section{Characterization of NWBFs over Local Field} 
In this section, we provide the characterization of NWBFs in $(FH^s(\Omega), FH^{-s}(\Omega))$. Firstly we need following two lemmas. 
\begin{lemma}\label{l 3.1} Given $s \in K$, let $\{T_{u(k)} \psi_0: k \in \mathbb N_0\}\cup \{T_{u(k)} \psi_\ell : k \in \mathbb N_0,  \; 1 \le \ell \le L\}$ be a Bessel sequence in $H^s(K)$. Then
\begin{eqnarray}
\nonumber && \sum_{k \in \mathbb N_0}\left |\langle  g, \psi_{0, 0, k}\rangle\right |^2+ \sum_{\ell =1}^L \sum_{j=0}^\infty \sum_{k \in \mathbb N_0}\left |\langle  g, \psi_{\ell, j, k}^s\rangle\right |^2\\\
\nonumber &&=\int_K \left| \hat g (\omega)\right|^2\left( \left| \hat \psi_0 (\omega)\right|^2+ \sum_{\ell =1}^L \sum_{j=0}^\infty  m^{-2js} \left| \hat \psi_\ell (\mathfrak p^{-j}\omega)\right|^2\right)d \omega\\\
\nonumber&&+\int_K \overline{\hat g (\omega)} \sum_{k \in \mathbb N_0}\hat g \left(\omega+u(k)\right)\\\
\label{3.1}&&\times\left( \hat \psi_0(\omega)\overline{\hat \psi_0 \left(\omega+u(k)\right)}+\sum_{\ell =1}^L \sum_{j=0}^k  m^{-2js} \hat \psi_\ell (\mathfrak p^{-j}\omega) \overline{\hat \psi_\ell \left(\omega+u(k)\right)}\right)
\end{eqnarray}
for $g \in \mathcal D$
\end{lemma}

{\it Proof.}  From Lemma 2.3, we have

\begin{eqnarray}
\nonumber && \sum_{k \in \mathbb N_0}\left |\langle  g, \psi_{0, 0, k}\rangle\right |^2+ \sum_{\ell =1}^L \sum_{j=0}^\infty \sum_{k \in \mathbb N_0}\left |\langle  g, \psi_{\ell, j, k}^s\rangle\right |^2\\\
\nonumber && =\int_{\mathfrak D}\left | \sum_{k \in \mathbb N_0}\hat g \left(\omega+u(k)\right)\overline{\hat \psi_0 \left(\omega+u(k)\right)}\right |^2d \omega\\\
\nonumber&&+\sum_{\ell =1}^L \sum_{j=0}^\infty m^{j(d-2s)}\int_{\mathfrak D}\left | \sum_{k \in \mathbb N_0}\hat g \left(\mathfrak p^j(\omega+u(k))\right)\overline{\hat \psi_\ell\left(\omega+u(k)\right)}\right |^2d \omega\\\
\nonumber && =\int_{\mathfrak D}\left (\sum_{k \in \mathbb N_0}\hat \psi_0 \left(\omega+u(k)\right)\overline{\hat g \left(\omega+u(k)\right)}\right )\left (\sum_{k \in \mathbb N_0}\hat g \left(\omega+u(k)\right)\overline{\hat \psi_0 \left(\omega+u(k)\right)}\right )d \omega\\\
\nonumber&&+\sum_{\ell =1}^L \sum_{j=0}^\infty m^{j(d-2s)}\int_{\mathfrak D}\left( \sum_{k \in \mathbb N_0}\hat \psi_\ell\left(\omega+u(k)\right)\overline{ \hat g \left(\mathfrak p^j(\omega+u(k))\right)}\right ) \\\
\nonumber && \times \left ( \sum_{k \in \mathbb N_0}\hat g \left(\mathfrak p^j(\omega+u(k))\right)\overline{\hat \psi_\ell\left(\omega+u(k)\right)}\right)\\\
\nonumber && =\int_{\mathfrak D}\left (\sum_{k \in \mathbb N_0}\hat \psi_0 \left(\omega+u(k)\right)\overline{\hat g \left(\omega+u(k)\right)}\right )E_0(\omega)d \omega\\\
\nonumber&&+\sum_{\ell =1}^L \sum_{j=0}^\infty m^{j(d-2s)}\int_{\mathfrak D}\left( \sum_{k \in \mathbb N_0}\hat \psi_\ell\left(\omega+u(k)\right)\overline{ \hat g \left(\mathfrak p^j(\omega+u(k))\right)}\right )E_{\ell, j} d \omega \\\
\label{3.2} && =R_1+R_2,
\end{eqnarray}
where $E_0(\cdot)= \sum_{k \in \mathbb N_0}\hat g \left(\omega+u(k)\right)\overline{\hat \psi_0 \left(\omega+u(k)\right)}$ and $E_{\ell, j}(\cdot)= \sum_{k \in \mathbb N_0}\hat g \left(\mathfrak p^j(\omega+u(k))\right)\overline{\hat \psi_\ell\left(\omega+u(k)\right)}$. Since $\{T_{u(k)} \psi_0: k \in \mathbb N_0\}$ is a Bessel sequence in $H^s(K)$ and $g \in \mathcal D$, it follows that $|E_0(\cdot)|\le [\hat g, \hat g]_{-s}^{\frac{1}{2}}(\cdot)[\hat \psi_0, \hat \psi_0]_{-s}^{\frac{1}{2}}(\cdot) < \infty$ by Lemma 2.2 and Lemma 2.5. Therefore 
$$\int_{\mathfrak D}\left | \sum_{k \in \mathbb N_0}\hat g \left(\omega+u(k)\right)\overline{\hat \psi_0 \left(\omega+u(k)\right)} E_0(\omega)\right | \le ||E_0||_{L^\infty(\mathfrak D)}\int_{\mathfrak D}[\hat g, \hat g]_{-s}^{\frac{1}{2}}(\omega)[\hat \psi_0, \hat \psi_0]_{-s}^{\frac{1}{2}}(\omega) < \infty$$,
hence by Fubini-Tonelli theorem
\begin{eqnarray*}
 &&\int_{\mathfrak D}\left (\sum_{k \in \mathbb N_0}\hat \psi_0 \left(\omega+u(k)\right)\overline{\hat g \left(\omega+u(k)\right)}\right )\left (\sum_{k \in \mathbb N_0}\hat g \left(\omega+u(k)\right)\overline{\hat \psi_0 \left(\omega+u(k)\right)}\right )d \omega\\\
&&=\int_{\mathfrak D}\hat \psi_0(\omega) \overline{\hat g(\omega )}\sum_{k \in \mathbb N_0}\hat g \left(\omega+u(k)\right)\overline{\hat \psi_0 \left(\omega+u(k)\right)}d \omega
\end{eqnarray*}
Moreover
\begin{eqnarray*}
&&\int_K |\hat \psi_0(\omega) \overline{\hat g(\omega )}|\sum_{k \in \mathbb N_0}|\hat g \left(\omega+u(k)\right)\overline{\hat \psi_0 \left(\omega+u(k)\right)}|d \omega\\\
 &&\le \int_{supp (\hat g)}\left(\sum_{k \in \mathbb N_0}|\hat g \left(\omega+u(k)\right)\overline{\hat \psi_0 \left(\omega+u(k)\right)} \right)^2 d \omega\\\
&&\le \int_{supp (\hat g)}[ \hat g, \hat g]_{-s}(\omega) [\hat \psi_0, \hat \psi_0]_s(\omega) d \omega\\\
&&< \infty
\end{eqnarray*}
since $[ \hat g, \hat g]_{-s}(\omega) [\hat \psi_0, \hat \psi_0]_s(\cdot) $ is essentially bounded by Lemma 2.2, we have
\begin{eqnarray}
\nonumber R_1&=&\int_K \hat \psi_0(\omega) \overline{\hat g(\omega )} \sum_{k \in \mathbb N_0}\hat g \left(\omega+u(k)\right)\overline{\hat \psi_0 \left(\omega+u(k)\right)}d \omega\\\
\label{3.3}&=&\int_K |\hat \psi_0(\omega)|^2| \hat  g(\omega )|^2 d \omega +\int_K \hat \psi_0(\omega) \overline{ g(\omega )} \sum_{k \in \mathbb N_0}\hat g \left(\omega+u(k)\right)\overline{\hat \psi_0 \left(\omega+u(k)\right)}d \omega
\end{eqnarray}
In order to complete the proof, we need to calculate $R_2$. Let us define $\tilde g$ by $\hat {\tilde g}(\cdot)=\hat g ((\mathfrak p)^{j}\cdot)$. Then  as $g \in \mathcal D$ and Lemma 2.5, we have
$$[\hat g ((\mathfrak p)^{j}\cdot), \hat g ((\mathfrak p)^{j}\cdot)](\cdot)\le C$$
Thus
$$|E_{i, j}(\cdot)|\le[\hat g ((\mathfrak p)^{j}\cdot), \hat g ((\mathfrak p)^{j}\cdot)]_{-s}^{\frac{1}{2}}(\cdot) [\hat \psi_\ell, \hat \psi_\ell ]_s^{\frac{1}{2}}< \infty.$$
Hence
$$R_2=\sum_{\ell=1}^L \sum{j=0}^\infty m^{j(d-2s)}\int_K \hat \psi_\ell(\omega) \overline{ g(\mathfrak p^j\omega )} \sum_{k \in \mathbb N_0}\hat g \left(\mathfrak p^j(\omega+u(k))\right)\overline{\hat \psi_\ell \left(\omega+u(k)\right)}d \omega
$$
Taking $A$ a bounded set in $K$ such that $supp (\hat g) \subset A$. So by Lemma 2.6, we get
$$A\cap (A+ \mathfrak p^j u(k))= \emptyset\quad \mbox{for} \; (j, k)\notin A_1 \times A_2\quad \mbox{with}\; k \ne 0,$$
where $A_1 \subset \mathbb N_0$ and $A_2 \subset \mathbb N_0 \setminus \{0\}$ are two finite sets. Therefore 
$$R_2=\sum_{\ell=1}^L \sum_{j\in A_1} m^{j(d-2s)}\int_K \hat \psi_\ell(\omega) \overline{ g(\mathfrak p^j\omega )} \sum_{k \in A_2}\hat g \left(\mathfrak p^j(\omega+u(k))\right)\overline{\hat \psi_\ell \left(\omega+u(k)\right)}d \omega $$
Denote $S=\bigcup_{k \in A_2 \cap \{0\}}(\bigcup_{j \in A_1}\left(\mathfrak p^j A+u(k)\right))$. Therefore, for each $(j, k)\in A_1 \times  A_2$, we have
\begin{eqnarray*}
&&\int_K \left |\hat \psi_\ell(\omega) \overline{ g(\mathfrak p^j\omega )} \hat g \left(\mathfrak p^j(\omega+u(k))\right)\overline{\hat \psi_\ell \left(\omega+u(k)\right)}\right |d \omega\\\
&& \le ||\hat g||^2_{L^\infty(K)}\int_{\mathfrak p^{-j}A}|\hat \psi_\ell(\omega) \hat \psi_\ell (\omega+u(k))| d\omega\\\
&& \le ||\hat g||^2_{L^\infty(K)}\left(\int_{\mathfrak p^{-j}A}|\hat \psi_\ell(\omega)|^2 d\omega\right)^{\frac{1}{2}}\left(\int_{\mathfrak p^{-j}A}|\hat \psi_\ell (\omega+u(k))|^2\right)^{\frac{1}{2}} d\omega\\\
&& \le ||\hat g||^2_{L^\infty(K)}\int_S |\hat \psi _\ell(\omega)|^2 d\omega
\end{eqnarray*}
Note that $1\le (\max_{\omega \in A}(1+|\omega|^2)^{-s})(1+|\omega|^2)^s$ for $\omega \in A$. Thus
\begin{eqnarray*}
&&\int_K \left |\hat \psi_\ell(\omega) \overline{ g(\mathfrak p^j\omega )} \hat g \left(\mathfrak p^j(\omega+u(k))\right)\overline{\hat \psi_\ell \left(\omega+u(k)\right)}\right |d \omega\\\
 &&\le \left(\max_{\omega \in A}(1+|\omega|^2)^{-s}\right)||\hat g||^2_{L^\infty(K)}\int_A |\hat \psi_\ell(\omega)|^2(1+|\omega|^2)^s d \omega\\\
 &&\le \left(\max_{\omega \in A}(1+|\omega|^2)^{-s}\right)||\hat g||^2_{L^\infty(K)}||\psi_\ell||^2_{H^s(K)}\\\
 &&<\infty
 \end{eqnarray*}
 On combining the formula given above, we get
\begin{eqnarray}
\nonumber R_2&=&\int_K\sum_{\ell=1}^L \sum_{j=0}^\infty m^{j(d-2s)}|\hat \psi_\ell(\omega)|^2|\hat g(\mathfrak p^j\omega)|^2 d\omega\\\
\nonumber &&+\int_K\sum_{\ell=1}^L \sum_{j=0}^\infty m^{j(d-2s)}\overline{\hat g(\mathfrak p^j\omega)}\hat \psi_\ell(\omega)\sum_{k \in \mathbb N_0} \hat g \left(\mathfrak p^j(\omega+u(k))\right)\overline{\hat \psi_\ell \left(\omega+u(k)\right)} d\omega\\\
\nonumber &=&\int_K\sum_{\ell=1}^L \sum_{j=0}^\infty m^{-2js}|\hat \psi_\ell(\mathfrak p^{-j}\omega)|^2|\hat g(\omega)|^2 d\omega\\\
\nonumber &&+\int_K\sum_{\ell=1}^L \sum_{j=0}^\infty m^{-2js}\overline{\hat g(\omega)}\hat \psi_\ell(\mathfrak p^{-j}\omega) \sum_{k \in \mathbb N_0} \hat g \left(\omega+\mathfrak p^{-j}u(k))\right)\overline{\hat \psi_\ell \left(\mathfrak p^{-j}\omega+u(k)\right)} d\omega\\\
\nonumber  &=&\int_K\sum_{\ell=1}^L \sum_{j=0}^\infty m^{-2js}|\hat \psi_\ell(\mathfrak p^{-j}\omega)|^2|\hat g(\omega)|^2 d\omega\\\
\label {3.4}&&+\int_K \overline{\hat g(\omega)}\sum_{k \in \mathbb N_0} \hat g \left(\omega+u(k)\right) \sum_{\ell=1}^L\sum_{j=0}^{\kappa(k)} m^{-2js}\hat \psi_\ell(\mathfrak p^{-j}\omega) \overline{\hat \psi_\ell \left(\mathfrak p^{-j}\omega+u(k)\right)} d\omega
\end{eqnarray}
using the definition of $\kappa(k)$. Hence using (\ref{3.2}), (\ref{3.3}) and (\ref{3.4}), we get (\ref{3.1}). This completes the proof of the lemma.
\begin{lemma}\label{l 3.2}
Given $s \in K$, let $X^s(\psi_0; \Psi)$ and $X^{-s}(\tilde \psi_0; \tilde \Psi)$ be a Bessel sequences in $H^s(K)$ and  $H^{-s}(K)$, respectively. Then for all $f, g \in \mathcal D$, we have 
\begin{eqnarray}
\nonumber &&\sum_{k \in \mathbb N_0} \langle f, \tilde \psi_{0,0,k}\rangle  \langle  \psi_{0,0,k}, g\rangle+\sum_{\ell=1}^L \sum_{j=0}^\infty\sum_{k \in \mathbb N_0}\langle f, \tilde \psi_{\ell,j,k}^{-s}\rangle  \langle  \psi_{\ell,j,k}^s, g\rangle\\\
\nonumber &&=\int_K \hat f (\omega)\overline{\hat g (\omega)}\left(\hat \psi_0 (\omega)\overline{ \hat {\tilde \psi}_0 (\omega)}+\sum_{\ell=1}^L \sum_{j=0}^\infty \hat \psi_\ell (\mathfrak p^{-j}\omega)\overline{ \hat {\tilde \psi}_\ell (\mathfrak p^{-j}\omega)}\right)d\omega\\\
\nonumber &&+\int_K \overline{\hat g (\omega)}\sum_{k \in \mathbb N_0}\hat f (\omega+u(k))\\\
\label{3.5}&&\times \left(\hat \psi_0 (\omega)\overline{ \hat {\tilde \psi}_0 (\omega+u(k))}+\sum_{\ell=1}^L \sum_{j=0}^{\kappa(k)} \hat \psi_\ell (\mathfrak p^{-j}\omega)\overline{ \hat {\tilde \psi}_\ell (\mathfrak p^{-j}(\omega+u(k))}\right)d\omega
\end{eqnarray}
\end{lemma}
{\it{proof}}: As $X^s(\psi_0; \Psi)$ and $X^{-s}(\tilde \psi_0; \tilde \Psi)$ are Bessel sequences in $H^s(K)$ and  $H^{-s}(K)$, respectively, the expression in (\ref{3.5}) is meaningful. Proceeding in a similar fashion as in Lemma \ref{l 3.1}, we have
\begin{eqnarray}
\nonumber &&\sum_{k \in \mathbb N_0} \langle f, \tilde \psi_{0,0,k}\rangle  \langle  \psi_{0,0,k}, g\rangle+\sum_{\ell=1}^L \sum_{j=0}^\infty\sum_{k \in \mathbb N_0}\langle f, \tilde \psi_{\ell,j,k}^{-s}\rangle  \langle  \psi_{\ell,j,k}^s, g\rangle\\\
\nonumber &&=\int_K \hat \psi_0(\omega)\overline{\hat g (\omega)}sum_{k \in \mathbb N_0}\hat f (\omega +u(k))\overline{\hat{ \tilde\psi}_0(\omega+u(k)) } d\omega\\\
\nonumber &&+\sum_{\ell=1}^L \sum_{j=0}^\infty q^j \int_K \hat \psi_\ell(\omega)\overline{\hat g (\mathfrak p^j\omega)}\sum_{k \in \mathbb N_0}\hat f (\mathfrak p^j(\omega +u(k)))\overline{\hat{ \tilde\psi}_\ell(\omega+u(k)) } d\omega\\\
\label{3.6}&& =I_1+I_2.
\end{eqnarray}
Note that
$$|\hat \psi_0(\cdot)\overline{\hat g (\cdot)}|\sum_{k \in \mathbb N_0}|\hat f (\cdot +u(k))\overline{\hat{ \tilde\psi}_0(\cdot+u(k)) }|\le [\hat f, \hat f]_{s}^{\frac{1}{2}}(\cdot)[\hat {\tilde \psi}_0, {\tilde \psi}_0]_{-s}^{\frac{1}{2}}(\cdot)[\hat g, \hat g]_{-s}^{\frac{1}{2}}(\cdot)[\hat { \psi}_0, { \psi}_0]_{s}^{\frac{1}{2}}(\cdot),$$
is bounded due to Lemma 2.2. Hence
\begin{eqnarray*}
&&\int_K |\hat \psi_0(\omega)\overline{\hat g (\omega)}|sum_{k \in \mathbb N_0}|\hat f (\omega +u(k))\overline{\hat{ \tilde\psi}_0(\omega+u(k)) } |d\omega\\\
&& \le \int_{supp (\hat g)} |\hat \psi_0(\omega)\overline{\hat g (\omega)}|sum_{k \in \mathbb N_0}|\hat f (\omega +u(k))\overline{\hat{ \tilde\psi}_0(\omega+u(k)) } |d\omega<\infty
\end{eqnarray*}
Therefore
\begin{equation}\label{3.7}
I_1=\int_K \hat f (\omega)\overline{\hat g (\omega)}\hat \psi_0 (\omega)\overline{ \hat {\tilde \psi}_0 (\omega)}+\int_K \hat \psi_0(\omega)\overline{\hat g (\omega)}sum_{k \in \mathbb N_0}\hat f (\omega +u(k))\overline{\hat{ \tilde\psi}_0(\omega+u(k)) } d\omega.
 \end{equation} 
 To complete the proof of the lemma, we need to discuss $I_2$. Let's break it into two parts, for $k=0$ and $k \ne 0$. Hence by  Cauchy-Schwartz inequality and Lemma 2.4,  we have
 \begin{eqnarray}
\nonumber &&\sum_{\ell=1}^L \sum_{j=0}^\infty |\hat \psi_\ell (\mathfrak p^{-j}\omega)\overline{ \hat {\tilde \psi}_\ell (\mathfrak p^{-j}\omega)} |\\\
\nonumber &&\le \left(\sum_{\ell=1}^L \sum_{j=0}^\infty  m^{-2js}|\hat \psi_\ell (\mathfrak p^{-j}\omega)|^2\right)^{\frac{1}{2}}  \left(\sum_{\ell=1}^L \sum_{j=0}^\infty  m^{2js}|\hat{\tilde  \psi}_\ell (\mathfrak p^{-j}\omega)|^2\right)^{\frac{1}{2}}\\\
\label{3.8}&&\le B_1 B_2.
 \end{eqnarray}
 Therefore
  \begin{eqnarray*}
&&\int_K |\hat f (\omega)\overline{\hat g (\omega)}|\sum_{\ell=1}^L \sum_{j=0}^\infty |\hat \psi_\ell (\mathfrak p^{-j}\omega)\overline{ \hat {\tilde \psi}_\ell (\mathfrak p^{-j}\omega)} |\\\
&&\le B_1B_2\left|supp (\hat f)\cap supp(\hat g)\right|\|\hat f\|_{L^\infty(K)}\|\hat g\|_{L^\infty(K)}\\\
&&<\infty
\end{eqnarray*}
 Fix a compact set $A \in K$ such that $supp (\hat f)\cap supp(\hat g) \subset A$. Using Lemma 2.6, it follows that
 \begin{equation}\label{3.9}
 A\cap (A+ \mathfrak p^j u(k))= \emptyset\quad \mbox{for} \; (j, k)\notin A_1 \times A_2\quad \mbox{with}\; k \ne 0
  \end{equation} 
 where $A_1 \subset \mathbb N_0$ and $A_2 \subset \mathbb N_0 \setminus \{0\}$ are two finite sets. With the same argument as applied to $R_2$, we have
\begin{eqnarray}
\nonumber && \int_K|\overline{\hat g (\mathfrak p^j\omega)}\hat f (\mathfrak p^j(\omega +u(k)))\hat \psi_\ell(\omega)\overline{\hat{ \tilde\psi}_\ell(\omega+u(k)) }| d\omega\\\
\nonumber &&\le \|\hat g\|_{L^\infty(K)}\|\hat f\|_{L^\infty(K)}  \left(\int_{\mathfrak p^{-j}A}|\hat \psi_\ell (\omega)|^2\right)^{\frac{1}{2}}  \left(\int_{\mathfrak p^{-j}A}|\hat{\tilde  \psi}_\ell (\omega+u(k))|^2\right)^{\frac{1}{2}}\\\
 \nonumber &&\le \|\hat g\|_{L^\infty(K)}\|\hat f\|_{L^\infty(K)}  \left(\int_{T}|\hat \psi_\ell (\omega)|^2\right)^{\frac{1}{2}}  \left(\int_{T}|\hat{\tilde  \psi}_\ell (\omega+u(k))|^2\right)^{\frac{1}{2}}\\\
 \nonumber &&\le \|\hat g\|_{L^\infty(K)}\|\hat f\|_{L^\infty(K)}  \left(\max_{\omega \in T}(1+|\omega|^2)\right)^{-s/2}\left(\max_{\omega \in T}(1+|\omega|^2)\right)^{s/2}\|\psi_\ell\|_{H^s(K)}\|\tilde {\psi}_\ell\|_{H^{-s}(K)}\\\
\label{3.10} &&<\infty
 \end{eqnarray}
 for $ (j, k)\in A_1 \times A_2$, where $T=\bigcup_{k\in A_2\cup \{0\}}(\bigcup_{j \in A_1}\mathfrak p^{-j}A+u(k))$. Using (\ref{3.8}) and (\ref{3.10}), we get
 \begin{eqnarray}
\nonumber I_2&=&\sum_{\ell=1}^L \sum_{j=0}^\infty \int_K\overline{\hat g (\omega)}\hat \psi_\ell(\mathfrak p^{-j}\omega)\sum_{k \in \mathbb N_0}\hat f (\omega +\mathfrak p^ju(k))\overline{\hat{ \tilde\psi}_\ell(\mathfrak p^{-j}(\omega+u(k))) } d\omega\\\
\nonumber &=& \int_K \sum_{\ell=1}^L \sum_{j=0}^\infty \overline{\hat g (\omega)}\hat \psi_\ell(\mathfrak p^{-j}\omega)\sum_{k \in \mathbb N_0}\hat f (\omega +\mathfrak p^ju(k))\overline{\hat{ \tilde\psi}_\ell(\mathfrak p^{-j}(\omega+u(k))) } d\omega\\\ 
\label{3.11}&=& \int_K \overline{\hat g (\omega)} \sum_{k \in \mathbb N_0}\hat f (\omega +u(k))\sum_{\ell=1}^L \sum_{j=0}^{\kappa(k)} \hat \psi_\ell(\mathfrak p^{-j}\omega)\overline{\hat{ \tilde\psi}_\ell(\mathfrak p^{-j}(\omega+u(k))) } d\omega\
\end{eqnarray}
On combining (\ref{3.6}), (\ref{3.7}) and (\ref{3.11}), we get (\ref{3.5}), which completes the proof of the lemma.

We now present a characterization of NWBFs in $(FH^s(\Omega), FH^{-s}(\Omega))$ in the form of the following theorem.
\begin{theorem}\label{t 3.1}
Given $s \in K$, let $FH^s(\Omega)$ and $FH^{-s}(\Omega)$ be reducing subspaces of $H^s(K)$ and $H^{-s}(K)$, respectively, $\psi_0 \in H^s(K), \; \tilde {\psi}_0 \in H^{-s}(K)$ and $\Psi \in H^s(K), \; \tilde {\Psi} \in H^{-s}(K)$. Suppose that $X^s(\psi_0, \Psi)$ and $X^{-s}(\tilde{\psi}_0, \tilde \Psi)$ are Bessel sequences in $FH^s(\Omega)$ and $FH^{-s}(\Omega)$, respectively. Then  $X^s(\psi_0, \Psi); X^{-s}(\tilde{\psi}_0, \tilde \Psi)$ is an NWBFs in $(FH^s(\Omega), FH^{-s}(\Omega))$ if and only if
\begin{equation}\label{3.12}
\hat \psi_0(\cdot)\overline{\hat{ \tilde\psi}_0(\cdot+u(k))}\sum_{\ell=1}^L \sum_{j=0}^{\kappa(k)} \hat \psi_\ell(\mathfrak p^{-j}\cdot)\overline{\hat{ \tilde\psi}_\ell(\mathfrak p^{-j}(\cdot+u(k))) } =\delta_{0, k}\quad a.e.\; on\; \Omega.
\end{equation}
\end{theorem}
{\it Proof}: As $\mathcal D \cap FH^s(\Omega)$ is dense in $FH^s(\Omega)$, then
$$X^s(\psi_0, \Psi); X^{-s}(\tilde{\psi}_0, \tilde \Psi)$$
ia an NWBF's in $(FH^s(\Omega), FH^{-s}(\Omega))$ iff for $f \in \mathcal D \cap FH^s(\Omega)$ and $g \in \mathcal D \cap FH^{-s}(\Omega)$ 

$$\sum_{k \in \mathbb N_0} \langle f, \tilde \psi_{0,0,k}\rangle  \langle  \psi_{0,0,k}, g\rangle+\sum_{\ell=1}^L \sum_{j=0}^\infty\sum_{k \in \mathbb N_0}\langle f, \tilde \psi_{\ell,j,k}^{-s}\rangle  \langle  \psi_{\ell,j,k}^s, g\rangle=\langle f, g \rangle$$
This is equivalent to
 \begin{eqnarray}
\nonumber &&\sum_{k \in \mathbb N_0} \langle (\hat f\chi_\Omega), \tilde \psi_{0,0,k}\rangle  \langle  \psi_{0,0,k}, (\hat g\chi_\Omega)\rangle+\sum_{\ell=1}^L \sum_{j=0}^\infty\sum_{k \in \mathbb N_0}\langle (\hat f\chi_\Omega), \tilde \psi_{\ell,j,k}^{-s}\rangle  \langle  \psi_{\ell,j,k}^s, (\hat g\chi_\Omega)\rangle\\\
\label{3.13}&&=\langle  (\hat f\chi_\Omega), (\hat g\chi_\Omega) \rangle
\end{eqnarray}
with $f, g \in \mathcal D$ as $\mathcal D \cap FH^s(\Omega)=\{(\hat h \chi_\Omega): h \in \mathcal D\}$. The expression of (\ref{3.13}) is well defined as $X^s(\psi_0, \Psi)$ and $X^{-s}(\tilde{\psi}_0, \tilde \Psi)$ are Bessel sequences in $H^s(K)$ and $H^{-s}(K)$. Hence using Lemma \ref{l 3.2}, we can write expression (\ref{3.13}) as
 \begin{eqnarray}
\nonumber &&\int_K \hat f (\omega)\overline{\hat g (\omega)}\chi_\Omega(\omega)\left(\hat \psi_0 (\omega)\overline{ \hat {\tilde \psi}_0 (\omega)}+\sum_{\ell=1}^L \sum_{j=0}^\infty \hat \psi_\ell (\mathfrak p^{-j}\omega)\overline{ \hat {\tilde \psi}_\ell (\mathfrak p^{-j}\omega)}\right)d\omega\\\
\nonumber &&+\int_K \overline{\hat g (\omega)}\chi_\Omega(\omega)\sum_{k \in \mathbb N_0}(\hat f \chi_\Omega)(\omega+u(k))\\\
\nonumber &&\times \left( \hat { \psi}_0 (\omega)\overline{\hat {\tilde \psi}_0 (\omega+u(k))}+\sum_{\ell=1}^L \sum_{j=0}^{\kappa(k)} \hat \psi_\ell (\mathfrak p^{-j}\omega)\overline{ \hat {\tilde \psi}_\ell (\mathfrak p^{-j}\omega+u(k))}\right)d\omega\\\
\label{3.14}&&=\int_K \hat f (\omega)\overline{\hat g (\omega)}\chi_\Omega(\omega)
\end{eqnarray}
with $ f, g \in \mathcal D$. Hence the expression (\ref{3.12}) leads (\ref{3.14}). It remains only to prove the converse statement. Suppose (\ref{3.14}) hold. Using Cauchy-Schwartz inequality, we have
\begin{eqnarray*}
&&| \hat { \psi}_0 (\cdot) \overline{\hat {\tilde \psi}_0 (\cdot+u(k))}|+\sum_{\ell=1}^L \sum_{j=0}^{\kappa(k)}| \hat \psi_\ell (\mathfrak p^{-j}\cdot)\overline{ \hat {\tilde \psi}_\ell (\mathfrak p^{-j}\cdot+u(k))}|\\\
&&\le \left(| \hat { \psi}_0 (\cdot) |^2+\sum_{\ell=1}^L \sum_{j=0}^{\infty}  m^{-2js}|\hat \psi_\ell (\mathfrak p^{-j}\cdot)|^2\right)^{{\frac{1}{2}}}\\\
&&+ \left(| \hat { \psi}_0 (\cdot+u(k)) |^2+\sum_{\ell=1}^L \sum_{j=0}^{\infty}  m^{-2js}|\hat \psi_\ell (\mathfrak p^{-j}(\cdot+u(k)))|^2\right)^{{\frac{1}{2}}}\\\
&& \le B_1 B_2 (1+|\cdot|^2)^{-s}(1+|\cdot+u(k)|^2)^{s}\\\
&&=C_k<\infty \end{eqnarray*}
for each $k \in \mathbb N_0$ using Lemma 2.4. Hence the series\\ $ \hat { \psi}_0 (\cdot) \overline{\hat {\tilde \psi}_0 (\cdot+u(k))}+\sum_{\ell=1}^L \sum_{j=0}^{\kappa(k)} \hat \psi_\ell (\mathfrak p^{-j}\cdot)\overline{ \hat {\tilde \psi}_\ell (\mathfrak p^{-j}\cdot+u(k))}$ converges absolutely a. e. on $K$ and is contained in $L^\infty(K)$. Therefore almost all points in $K$ are its Lebesgue points. Next we consider two cases. When $k=0$. Let $\omega_0 \ne 0$ be a Lebesgue point of 
$ \hat { \psi}_0 (\cdot) \overline{\hat {\tilde \psi}_0 (\cdot)}+\sum_{\ell=1}^L \sum_{j=0}^{\infty} \hat \psi_\ell (\mathfrak p^{-j}\cdot)\overline{ \hat {\tilde \psi}_\ell (\mathfrak p^{-j}\cdot)}$ and $\chi_\Omega(\cdot)$. Fix $f$ and $g $ for $0< \epsilon<u(1)$, we have
$$\hat f(\cdot)=\hat g(\cdot)=\frac{\chi_{B(\omega_0, \epsilon)}}{\sqrt{|B(\omega_0, \epsilon)|}}$$
in (\ref{3.14}), where $B(\omega_0, \epsilon)$ is an open ball centred at $\omega_0$ and radius $\epsilon$. Therefore
\begin{eqnarray*}
&&\frac{1}{|B(\omega_0, \epsilon)|}\int_{B(\omega_0, \epsilon)}\chi_{\Omega}(\omega)\left( \hat { \psi}_0 (\omega) \overline{\hat {\tilde \psi}_0 (\omega)}+\sum_{\ell=1}^L \sum_{j=0}^{\infty} \hat \psi_\ell (\mathfrak p^{-j}\omega)\overline{ \hat {\tilde \psi}_\ell (\mathfrak p^{-j}\omega)}\right)d \omega\\\
&&=\frac{1}{|B(\omega_0, \epsilon)|}\int_{B(\omega_0, \epsilon)}\chi_{\Omega}(\omega) d\omega
\end{eqnarray*}
letting $\epsilon \to 0$, we have
$$ \hat { \psi}_0 (\omega_0) \overline{\hat {\tilde \psi}_0 (\omega_0)}+\sum_{\ell=1}^L \sum_{j=0}^{\infty} \hat \psi_\ell (\mathfrak p^{-j}\omega_0)\overline{ \hat {\tilde \psi}_\ell (\mathfrak p^{-j}\omega_0)}$$
For $ k \ne 0$, we fix $k_0\in \mathbb N_0$ and take $f$ and $g$
$$\hat f(\cdot+u(k))=\hat g(\cdot)=\frac{\chi_{B(\omega_0, \epsilon)}}{\sqrt{|B(\omega_0, \epsilon)|}}$$
in (\ref{3.14}), with $0< \epsilon <\frac{1}{2}$. Therefore
\begin{eqnarray*}
&&\frac{1}{|B(\omega_0, \epsilon)|}\int_{B(\omega_0, \epsilon)}\chi_{\Omega}(\omega)\\\
&&\times \left( \hat { \psi}_0 (\omega) \overline{\hat {\tilde \psi}_0 (\omega+u(k))}+\sum_{\ell=1}^L \sum_{j=0}^{\kappa(k)} \hat \psi_\ell (\mathfrak p^{-j}\omega)\overline{ \hat {\tilde \psi}_\ell (\mathfrak p^{-j}\omega+u(k))}\right)d \omega=0
\end{eqnarray*}
letting $\epsilon \to 0$ and using Lebesgue differentiation theorem, we get
$$ \hat { \psi}_0 (\omega_0+u(k_0)) \overline{\hat {\tilde \psi}_0 (\omega_0)}+\sum_{\ell=1}^L \sum_{j=0}^{\infty} \hat \psi_\ell (\mathfrak p^{-j}\omega_0)\overline{ \hat {\tilde \psi}_\ell (\mathfrak p^{-j}\omega_0+u(k_0))}$$
Hence we obtain (\ref{3.12})  by using the arbitrariness of $\omega_0$ and $k_0$, which completes the proof of the theorem.

{\small

\end{document}